\documentstyle[a4,txmac,%
case,%
mypic,%
epsf,%
twoside,%
nocaphead,%
amssymb,times,mathptm]{article}

\makeatletter

\advance\oddsidemargin by -1.7cm
\advance\evensidemargin by -1.7cm
\advance\textwidth by 3.4cm

\def\mynewtheo#1#2{%
\newtheorem{@#1}{#2}[section]%
\newenvironment{#1}{\begin{@#1}\rm}{\end{@#1}}}
 
\mynewtheo{lemma}{Lemma}
\mynewtheo{exer}{Exercise}
\mynewtheo{theo}{Theorem}
\mynewtheo{rem}{Remark}
\mynewtheo{defi}{Definition}
\mynewtheo{conj}{Conjecture}
\mynewtheo{corr}{Corollary}
\mynewtheo{prop}{Proposition}
\mynewtheo{question}{Question}
\mynewtheo{exam}{Example}

\def\Mvariable#1{\mathord {\operator@font #1}}
\def\Muserfunction#1{\mathord { #1}}
\def\inx{\mathop {\operator@font ind}\mathord{}}
\def\mpb{\mathop {\operator@font mpb}\mathord{}}
\def\spn{\mathop {\operator@font span}\mathord{}}
\def\len{\mathop {\operator@font len}\mathord{ }}

\newenvironment{theorem}{\begin{theo}}{\end{theo}}

\def\eqref#1{\mbox{(\protect\reference{#1}})}
\def\proof{\@ifnextchar[{\@proof}{\@proof[\unskip]}}
\def\@proof[#1]{\noindent{\bf Proof #1.}\enspace}

\begin{document}

\author{A. Stoimenow\\[2mm]
\small Research Institute for Mathematical Sciences, \\
\small Kyoto University, Kyoto 606-8502, Japan\\
\small e-mail: {\tt stoimeno@kurims.kyoto-u.ac.jp}\\
\small WWW: {\hbox{\web|http://www.kurims.kyoto-u.ac.jp/~stoimeno/|}}
}

\title{\large\bf \uppercase{5-moves and Montesinos links}
\\[4mm] {\small\it This is a preprint.
I would be grateful for any comments and corrections.
}}

% On values of the Jones polynomial and genus-minimizing diagrams

\date{\large Current version: \today\ \ \ First version:
\makedate{14}{3}{2006}}

\maketitle

\def\chrd#1#2{\picline{1 #1 polar}{1 #2 polar}}
\def\arrow#1#2{\picvecline{1 #1 polar}{1 #2 polar}}

\def\labch#1#2#3{\chrd{#1}{#2}\picputtext{1.3 #2 polar}{$#3$}}
\def\labar#1#2#3{\arrow{#1}{#2}\picputtext{1.3 #2 polar}{$#3$}}
\def\labbr#1#2#3{\arrow{#1}{#2}\picputtext{1.3 #1 polar}{$#3$}}

\def\rottab#1#2{%
\expandafter\advance\csname c@table\endcsname by -1\relax
\centerline{%
\rbox{\centerline{\vbox{\setbox1=\hbox{#1}%
\centerline{\mbox{\hbox to \wd1{\hfill\mbox{\vbox{{%
\caption{#2}}}}\hfill}}}%
\vskip9mm
\centerline{
\mbox{\copy1}}}}%
}%
}%
}

\def\GD{\szCD{6mm}}
\def\szCD#1#2{{\let\@nomath\@gobble\small\diag{#1}{2.4}{2.4}{
  \picveclength{0.27}\picvecwidth{0.1}
  \pictranslate{1.2 1.2}{
    \piccircle{0 0}{1}{}
    #2
}}}}
\def\CD{\szCD{4mm}}

\def\point#1{{\picfillgraycol{0}\picfilledcircle{#1}{0.08}{}}}
\def\labpt#1#2#3{\pictranslate{#1}{\point{0 0}\picputtext{#2}{$#3$}}}
\def\vrt#1{{\picfillgraycol{0}\picfilledcircle{#1}{0.09}{}}}

\def\@dcont{}
\def\svCD#1{\ea\glet\csname #1\endcsname\@dcont}
\def\rsCD#1{\ea\glet\ea\@dcont\csname #1\endcsname\ea\glet
\csname #1\endcsname\relax}

\def\addCD#1{\ea\gdef\ea\@dcont\ea{\@dcont #1}}
\def\drawCD#1{\szCD{#1}{\@dcont}}

\def\noloop{{\diag{0.5cm}{0.5}{1}{\picline{0.25 0}{0.25 1}}}}

\def\vrt#1{{\picfillgraycol{0}\picfilledcircle{#1}{0.09}{}}}

\def\ReidI#1#2{
  \diag{0.5cm}{0.9}{1}{
    \pictranslate{0.4 0.5}{
      \picscale{#11 #21}{
        \picmultigraphics[S]{2}{1 -1}{
	  \picmulticurve{-6 1 -1 0}{0.5 -0.5}{0.5 0}{0.1 0.3}{-0.2 0.3}
	} %picmultigraphics
	\piccirclearc{-0.2 0}{0.3}{90 270}
      }%picscale
    }%pictranslate
  }%diag
}
\def\Pos#1#2{{\diag{#1}{1}{1}{#2
\picmultiline{-5 1 -1 0}{0 1}{1 0}
\picmultiline{-5 1 -1 0}{0 0}{1 1}
}}}
\def\Neg#1#2{{\diag{#1}{1}{1}{#2
\picmultiline{-5 1 -1 0}{0 0}{1 1}
\picmultiline{-5 1 -1 0}{0 1}{1 0}
}}}
\def\Nul#1#2{{\diag{#1}{1}{1}{#2
\piccirclearc{0.5 1.4}{0.7}{-135 -45}
\piccirclearc{0.5 -0.4}{0.7}{45 135}
}}}
\def\Inf#1#2{{\diag{#1}{1}{1}{#2
\piccirclearc{0.5 1.4 x}{0.7}{135 -135}
\piccirclearc{0.5 -0.4 x}{0.7}{-45 45}
}}}
\def\pos{\Pos{0.5em}{\piclinewidth{10}}}
\def\neg{\Neg{0.5em}{\piclinewidth{10}}}
\def\nul{\Nul{0.5em}{\piclinewidth{10}}}
\def\iinf{\Inf{0.5em}{\piclinewidth{10}}}

\def\pt#1{{\picfillgraycol{0}\picfilledcircle{#1}{0.1}{}}}
\def\ppt#1{{\picfillgraycol{0}\picfilledcircle{#1}{0.06}{}}}

%usage ODD # of  coords, then {}
\def\@curvepath#1#2#3{%
  \@ifempty{#2}{\piccurveto{#1 }{@stc}{@std}#3}%
    {\piccurveto{#1 }{#2 }{#2  #3  0.5 conv}
    \@curvepath{#3}}%
}
\def\curvepath#1#2#3{%
  \piccurve{#1 }{#2 }{#2 }{#2  #3  0.5 conv}%
  \picPSgraphics{/@stc [ #1  #2  -1 conv ] $ D /@std [ #1  ] $ D }%
  \@curvepath{#3}%
}

%usage EVEN # of  coords, then {}
\def\@opencurvepath#1#2#3{%
  \@ifempty{#3}{\piccurveto{#1 }{#1 }{#2 }}%
    {\piccurveto{#1 }{#2 }{#2  #3  0.5 conv}\@opencurvepath{#3}}%
}
\def\opencurvepath#1#2#3{%
  \piccurve{#1 }{#2 }{#2 }{#2  #3  0.5 conv}%
  \@opencurvepath{#3}%
}

\long\def\@makecaption#1#2{%
   % \tm
   \vskip 10pt
   {\let\label\@gobble
   \let\ignorespaces\@empty
   \xdef\@tempt{#2}%
   %\typeout{`#2'}%
   }%
   \ea\@ifempty\ea{\@tempt}{%
   \setbox\@tempboxa\hbox{%
      \fignr#1#2}%
      }{%
   \setbox\@tempboxa\hbox{%
      {\fignr#1:}\capt\ #2}%
      }%
   \ifdim \wd\@tempboxa >\captionwidth {%
      \rightskip=\@captionmargin\leftskip=\@captionmargin
      \unhbox\@tempboxa\par}%
   \else
      \hbox to\captionwidth{\hfil\box\@tempboxa\hfil}%
   \fi}%
\def\fignr{\small\sffamily\bfseries}%
\def\capt{\small\sffamily}%

\newdimen\@captionmargin\@captionmargin2\parindent\relax
\newdimen\captionwidth\captionwidth\hsize\relax

\def\nin{\not\in}
\def\bQ{{\Bbb Q}}
\def\bC{{\Bbb C}}
\def\bR{{\Bbb R}}
\def\bN{{\Bbb N}}
\def\bZ{{\Bbb Z}}
\def\cE{{\cal E}}
\def\cK{{\cal K}}
\def\cI{{\cal I}}
\def\cR{{\cal R}}
\def\cV{{\cal V}}
\def\cO{{\cal O}}
\def\hD{{\hat D}}
\def\cD{{\cal D}}
\def\cX{{\cal X}}
\def\hK{{\hat K}}
\def\bm{\bar t'_2}
\let\bt\bm 

\def\pr{\text{\rm pr}\,}
\def\ncap{\not\mathrel{\cap}}
\def\|{\mathrel{\kern1.5pt\Vert\kern1.5pt}}
\def\lra{\longrightarrow}
\let\ds\displaystyle
\let\reference\ref
\def\so{\Rightarrow}
\let\x\exists
\let\fa\forall
\let\es\enspace
\def\llra{\longleftrightarrow}

\let\x\exists
\let\sg\sigma
\let\tl\tilde
\let\ap\alpha
\let\dl\delta
\let\Dl\Delta
\let\be\beta
\let\gm\gamma
\let\Gm\Gamma
\let\nb\nabla
\let\Lm\Lambda
\let\sm\setminus
\let\dt\det

\let\eps\varepsilon
\let\ul\underline
\let\ol\overline
\def\md{\min\deg}
\def\Md{\max\deg}
\def\Mcf{\max{\operator@font cf}}
\let\Mc\Mcf
\def\sgn{{\operator@font sgn}}
\def\TM{$^\text{\raisebox{-0.2em}{${}^\text{TM}$}}$}
\newcommand{\wt}{\widetilde}

\def\ssim{\stackrel{\ds \sim}{\vbox{\vskip-0.2em\hbox{$\scriptstyle *$}}}}

\def\bysame{\same[\kern2cm]\,}
\def\qed{\hfill\@mt{\Box}}
\def\@mt#1{\ifmmode#1\else$#1$\fi}

\def\proof{\@ifnextchar[{\@proof}{\@proof[\unskip]}}
\def\@proof[#1]{\noindent{\bf Proof #1.}\enspace}

\def\myfrac#1#2{\raisebox{0.2em}{\small$#1$}\!/\!\raisebox{-0.2em}{\small$#2$}}
\newcommand{\mybr}[2]{\text{$\Bigl\lfloor\mbox{%
\small$\displaystyle\frac{#1}{#2}$}\Bigr\rfloor$}}
\def\mybrtwo#1{\mbox{\mybr{#1}{2}}}

\def\epsfs#1#2{{\ifautoepsf\unitxsize#1\relax\else
\epsfxsize#1\relax\fi\epsffile{#2.eps}}}
\def\epsfsv#1#2{{\vcbox{\epsfs{#1}{#2}}}}
\def\vcbox#1{\setbox\@tempboxa=\hbox{#1}\parbox{\wd\@tempboxa}{\box
  \@tempboxa}}

\def\eepsfs{\@ifnextchar[{\@eepsfs}{\@@eepsfs}}
\def\@eepsfs[#1]#2{\uu{\ifautoepsf\unitxsize#1\relax\else
\epsfxsize#1\relax\fi\epsffile{#2.eps}}}
\def\@@eepsfs#1{\uu{\epsffile{#1.eps}}}

\def\uu#1{\setbox\@tempboxa=\hbox{#1}\@tempdima=-0.5\ht\@tempboxa
\advance\@tempdima by 0.5em\raise\@tempdima\box\@tempboxa}

\def\@test#1#2#3#4{%
  \let\@tempa\go@
  \@tempdima#1\relax\@tempdimb#3\@tempdima\relax\@tempdima#4\unitxsize\relax
  \ifdim \@tempdimb>\z@\relax
    \ifdim \@tempdimb<#2%
      \def\@tempa{\@test{#1}{#2}}%
    \fi
  \fi
  \@tempa
}

\def\go@#1\@end{}
\newdimen\unitxsize
\newif\ifautoepsf\autoepsftrue

\unitxsize4cm\relax
\def\epsfsize#1#2{\epsfxsize\relax\ifautoepsf
  {\@test{#1}{#2}{0.1 }{4   }
		{0.2 }{3   }
		{0.3 }{2   }
		{0.4 }{1.7 }
		{0.5 }{1.5 }
		{0.6 }{1.4 }
		{0.7 }{1.3 }
		{0.8 }{1.2 }
		{0.9 }{1.1 }
		{1.1 }{1.  }
		{1.2 }{0.9 }
		{1.4 }{0.8 }
		{1.6 }{0.75}
		{2.  }{0.7 }
		{2.25}{0.6 }
		{3   }{0.55}
		{5   }{0.5 }
		{10  }{0.33}
		{-1  }{0.25}\@end
		\ea}\ea\epsfxsize\the\@tempdima\relax
		\fi
		}

% typing `equation' breaks my fingers!
\newenvironment{eqn}{\begin{equation}}{\end{equation}\@ignoretrue}

% eqlabel=(#1)
\newenvironment{myeqn*}[1]{\begingroup\def\@eqnnum{\reset@font\rm#1}%
\xdef\@tempk{\arabic{equation}}\begin{equation}\edef\@currentlabel{#1}}
{\end{equation}\endgroup\setcounter{equation}{\@tempk}\ignorespaces}

% eqlabel=#1
\newenvironment{myeqn}[1]{\begingroup\let\eq@num\@eqnnum
\def\@eqnnum{\bgroup\let\r@fn\normalcolor % an extremely UGLY hack !!!
\def\normalcolor####1(####2){\r@fn####1#1}%
%\show\reset@font
\eq@num\egroup}%
\xdef\@tempk{\arabic{equation}}\begin{equation}\edef\@currentlabel{#1}}
{\end{equation}\endgroup\setcounter{equation}{\@tempk}\ignorespaces}

% eqlabel=(eqnnr) \qed
\newenvironment{myeqn**}{\begin{myeqn}{(%\arabic{equation}
%\show\theequation
\theequation)\es\es\mbox{\qed}}\edef\@currentlabel{\theequation}}
{\end{myeqn}\stepcounter{equation}}

\let\diagram\diag

\def\boxed#1{\diagram{1em}{1}{1}{\picbox{0.5 0.5}{1.0 1.0}{#1}}}

\def\rato#1{\hbox to #1{\rightarrowfill}}
\def\arrowname#1{{\enspace
\setbox7=\hbox{F}\setbox6=\hbox{%
\setbox0=\hbox{\footnotesize $#1$}\setbox1=\hbox{$\to$}%
\dimen@\wd0\advance\dimen@ by 0.66\wd1\relax
$\stackrel{\rato{\dimen@}}{\copy0}$}%
\ifdim\ht6>\ht7\dimen@\ht7\advance\dimen@ by -\ht6\else
\dimen@\z@\fi\raise\dimen@\box6\enspace}}

\def\contr{\diagram{1em}{0.6}{1}{\piclinewidth{35}%
\picstroke{\picline{0.5 1}{0.2 0.4}%
\piclineto{0.6 0.6}\picveclineto{0.3 0}}}}

\def\abstractname{}

\parskip5pt plus 1pt minus 2pt
\parindent\z@

{\let\@noitemerr\relax
\vskip-2.7em\kern0pt\begin{abstract}
\noindent{\bf Abstract.}\enspace
We classify the Montesinos links up to mutation and 5-move equivalence,
and obtain from this a Jones and Kauffman polynomial test for a
Montesinos link.
\\[1mm]
{\it Keywords:} rational link, Montesinos link, pretzel link, $n$-move,
Jones polynomial, Kauffman polynomial\\
{\it AMS subject classification:} 57M25 (primary)
\end{abstract}
}
\vspace{1mm}

{\parskip0.2mm\tableofcontents}
\vspace{7mm}

\section{Introduction}

$k$-moves are a family of natural operations on knot and link diagrams.
They were introduced by Nakanishi, and were studied first more
systematically by Przytycki \cite{Prz}. The following picture shows
as an example the move for $k=3$:
\begin{eqn}\label{3m}
\diag{4mm}{2}{2}{
  \piccirclearc{1 -1}{1.41}{45 135}
  \piccirclearc{1 3}{1.41}{-135 -45}
} 
\quad\llra\quad
\diag{6mm}{3.5}{1}{
  \piclinewidth{40}
  \picmultigraphics{3}{1 0}{
    \picline{0.2 0.8}{0.8 0.2}
    \picmultiline{-8.5 1 -1.0 0}{0.2 0.2}{0.8 0.8}
  }
  \picmultigraphics{2}{1 0}{
    \piccirclearc{1 0.6}{0.28}{45 135}
    \piccirclearc{1 0.4}{0.28}{-135 -45}
  }
  \picline{-0.2 -0.2}{0.3 0.3}
  \picline{-0.2 1.2}{0.3 0.7}
  \picline{3.2 -0.2}{2.7 0.3}
  \picline{3.2 1.2}{2.7 0.7}
}\,.
\end{eqn}

We call two links \em{$k$-move equivalent}, or simpler
\em{$k$-equivalent}, if they are related by a sequence of
$k$-moves (or their inverses). The case $k=2$ is the usual
crossing change. The cases $k=3,4$ were connected to two
long-standing problems. The $3$-move conjecture stated that
all links are $3$-move equivalent to a trivial link. It
was refuted only recently in \cite{DP}. The conjecture is
known to be true for many links, for example 3-algebraic
links \cite{PT}, links of braid index 4 and 5 \cite{Chen} (latter
with the exception of one equivalence class, which later provided
a counterexample), and knots of weak genus two \cite{gen2}.
The $4$-move conjecture states that all knots are $4$-move equivalent
to the unknot. This conjecture remains open, though partial
confirmations exist, and (in the general case) counterexamples
are suspected (see \cite{Ask,Prz2}, and also \cite{gen2}).

The difficulty with the cases $k=3,4$ is that, while equivalence
classes are rather few and large, (possibly exactly because of this)
it is impossible to obtain essential information on them using
polynomial (or other easy to compute) invariants.

However, the case $k=5$ is different. Now there are some invariants
that come from the Jones $V$ and Kauffman $F$ polynomial. Unfortunately,
the equivalence classes become many, and their intersection with a
meaningful family of links is difficult to describe. Recently,
Przytycki-D\c abkowski-Ishiwata (see \cite{Ishiwata})
succeeded in determining the classes for rational (or 2-bridge)
links. The aim of this paper is to extend their result to
Montesinos links (up to mutation); see theorems \reference{thq}
and \reference{thq2}.

The motivation for our work is to gain a condition on the polynomial
invariants of Montesinos links. A few such conditions were obtained for
partial classes, but all conditions for a general Montesinos link we
have so far are inefficient, or not easy to test.

In \cite{LickThis} semiadequate links were introduced. It was observed
that Montesinos links are semiadequate, and that for such links one
of the leading or trailing coefficients of the Jones polynomial $V$
must be $\pm 1$. (In \cite{ntriv} we understood also coefficients
2 and 3, and in particular proved that the Jones polynomial is
non-trivial.) However, this property is not helpful as a Montesinos
link test, because it is satisfied for many (other) links.
Similarly impractical is the semiadequacy condition on the Kauffman
$F$ polynomial of \cite{Thistle}. (The simplest knots whose polynomial
shows negative ``critical line'' coefficients on either side have
15 crossings; see \cite{ntriv}.) For the Alexander polynomial no
conditions whatsoever are known, that is, it is possible that every
admissible Alexander polynomial is realized by a Montesinos link.
In \cite{LickThis} it was shown how to determine the crossing number
of a Montesinos link from the $V$ and $F$ polynomial, but a(n
extensive, though systematic) diagram verification is not promising
either as a Montesinos link test. 

Now, in contrast, since we can evaluate the $5$-move invariants of
the Jones and Kauffman polynomial on the classes of Montesinos links
we obtain (theorems \reference{thfr} and \reference{thfs}), we gain
a condition on these polynomials, which turns out easy to
verify. (It makes no assumption on the diagram we perform it on,
except, of course, that one can evaluate the polynomials.)
We will give some examples that show how to apply our test.

\section{Preliminaries}

\subsection{Polynomial invariants}

Let $V$ be the Jones polynomial \cite{Jones}, $Q$ the BLMH
polynomial and $F$ the Kauffman polynomial \cite{Kauffman}.
We give a basic description of these invariants.

\begin{figure}[htb]
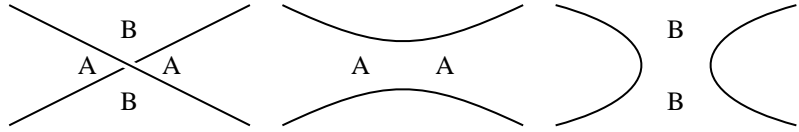

\[
\diag{8mm}{4}{2}{
   \picline{0 0}{4 2}
   \picmultiline{-5.0 1 -1.0 0}{0 2}{4 0}
   \picputtext{2.7 1}{A}
   \picputtext{1.3 1}{A}
   \picputtext{2 1.6}{B}
   \picputtext{2 0.4}{B}
} \quad
\diag{8mm}{4}{2}{
   \pictranslate{2 1}{
       \picmultigraphics[S]{2}{1 -1}{
           \piccurve{-2 1}{-0.3 0.2}{0.3 0.2}{2 1}
       }
   }
   \picputtext{2.7 1}{A}
   \picputtext{1.3 1}{A}
} \quad
\diag{8mm}{4}{2}{
   \pictranslate{2 1}{
       \picmultigraphics[S]{2}{-1 1}{
           \piccurve{2 -1}{0.1 -0.5}{0.1 0.5}{2 1}
       }
   }
   \picputtext{2 1.6}{B}
   \picputtext{2 0.4}{B}
}
\]
\caption{\label{figsplit}The A- and B-corners of a
crossing, and its both splittings. The corner A (resp. B)
is the one passed by the overcrossing strand when rotated 
counterclockwise (resp. clockwise) towards the undercrossing 
strand. A type A (resp.\ B) splitting is obtained by connecting 
the A (resp.\ B) corners of the crossing.}
\end{figure}

Recall first the construction of the Kauffman bracket in
\cite{Kauffman2}. The Kauffman bracket $[D]$ of a
diagram $D$ is a Laurent polynomial in a variable
$A$, obtained by summing over all states the terms
\begin{eqn}\label{eq_12}
A^{\#A-\#B}\,\left(-A^2-A^{-2}\right)^{|S|-1}\,.
\end{eqn}
A state is a choice of splittings of type $A$ or 
$B$ for any single crossing (see figure \ref{figsplit}), 
$\#A$ and $\#B$ denote the number of
type A (resp. type B) splittings and $|S|$ the
number of (disjoint) circles obtained after all
splittings in a state.

The Jones polynomial of a link $L$ is related to the
Kauffman bracket of some diagram $D$ of $L$ by
\begin{eqn}\label{conv}
V_L(t)\,=\,\left(-t^{-3/4}\right)^{-w(D)}\,[D]
\raisebox{-0.6em}{$\Big |_{A=t^{-1/4}}$}\,.
\end{eqn}

The \em{Kauffman polynomial} \cite{Kauffman} $F$ is usually defined
via a regular isotopy invariant $\Lm(a,z)$ of unoriented links.

We use here a slightly different convention for the variables
in $F$, differing from \cite{Kauffman,Thistle} by the interchange
of $a$ and $a^{-1}$. Thus in particular we have the relation
$F(D)(a,z)=a^{w(D)}\Lm(D)(a,z)$, where $w(D)$ is the writhe
of a link diagram $D$, and $\Lm(D)$ is the writhe-unnormalized
version of $F$. $\Lm$ is given in our convention by the properties
\begin{eqn}\label{(2)}
\begin{array}{c}
\Lambda\bigl(\Pos{0.5cm}{}\bigr)\ +\ \Lambda\bigl(\Neg{0.5cm}{}\bigr)\ =\ z\ \bigl(\ 
\Lambda\bigl(\Nul{0.5cm}{}\bigr)\ +\ \Lambda\bigl(\Inf{0.5cm}{}\bigr)\ \bigr)\,,\\[2mm]
\Lambda\bigl(\ \ReidI{-}{-}\bigr) = a^{-1}\ \Lambda\bigl(\noloop\bigr);\quad
\Lambda\bigl(\ \ReidI{-}{ }\bigr) = a\ \Lambda\bigl(\noloop\bigr)\,,\\[2mm]
\Lambda\bigl(\,\mbox{\Large $\bigcirc$}\,\bigr) = 1\,.
\end{array}
\end{eqn}

The BLMH polynomial $Q$ is most easily specified by $Q(z)=F(1,z)$.

\subsection{Families of links}

In the following we define rational, pretzel and Montesinos
links according to Conway \cite{Conway}.

\begin{figure}[htb]
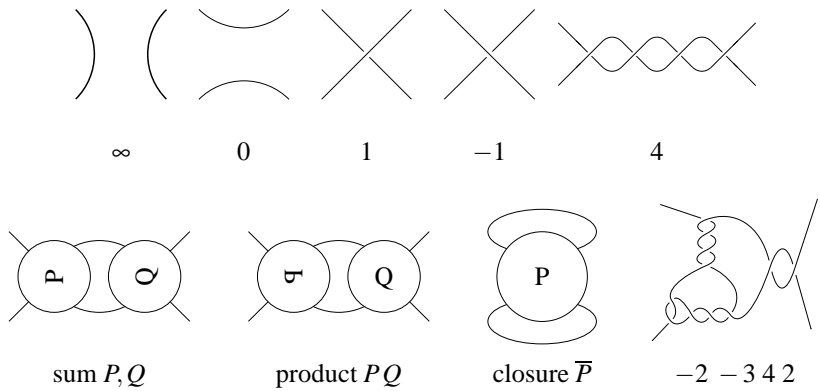

\[
\begin{array}{*5c}
\diag{6mm}{2}{2}{
  \piccirclearc{-1 1}{1.41}{-45 45}
  \piccirclearc{3 1}{1.41}{135 -135}
} & 
\diag{6mm}{2}{2}{
  \piclinewidth{50}
  \piccirclearc{1 -1}{1.41}{45 135}
  \piccirclearc{1 3}{1.41}{-135 -45}
} &
\diag{6mm}{2}{2}{
  \piclinewidth{50}
  \picmultiline{-9 1 -1.0 0}{2 0}{0 2}
  \picmultiline{-9 1 -1.0 0}{0 0}{2 2}
} &
\diag{6mm}{2}{2}{
  \piclinewidth{50}
  \picmultiline{-9 1 -1.0 0}{0 0}{2 2}
  \picmultiline{-9 1 -1.0 0}{2 0}{0 2}
} &
\diag{6mm}{4}{1}{
  \piclinewidth{50}
  \picmultigraphics{4}{1 0}{
    \picline{0.2 0.8}{0.8 0.2}
    \picmultiline{-6 1 -1.0 0}{0.2 0.2}{0.8 0.8}
  }
  \picmultigraphics{3}{1 0}{
    \piccirclearc{1 0.6}{0.28}{45 135}
    \piccirclearc{1 0.4}{0.28}{-135 -45}
  }
  \picline{-0.2 -0.2}{0.3 0.3}
  \picline{-0.2 1.2}{0.3 0.7}
  \picline{4.2 -0.2}{3.7 0.3}
  \picline{4.2 1.2}{3.7 0.7}
} \\[9mm]
\infty & 0 & 1 & -1 & 4 
\end{array}
\]
\[
\begin{array}{c@{\qquad}c@{\qquad}c@{\qquad}c}
\diag{6mm}{4}{2}{
  \piclinewidth{50}
  \picline{0 0}{0.5 0.5}
  \picline{0 2}{0.5 1.5}
  \picline{4 0}{3.5 0.5}
  \picline{4 2}{3.5 1.5}
  \piccirclearc{2 0.5}{1.3}{45 135}
  \piccirclearc{2 1.5}{1.3}{-135 -45}
  \pictranslate{1 1}{
    \picrotate{90}{
      \picscale{1 1}{
	\picfilledcircle{0 0}{0.8}{P}
      }
    }
  }
  \pictranslate{3 1}{
    \picrotate{90}{
      \picscale{1 1}{
	\picfilledcircle{0 0}{0.8}{Q}
      }
    }
  }
} &
\diag{6mm}{4}{2}{
  \piclinewidth{50}
  \picline{0 0}{0.5 0.5}
  \picline{0 2}{0.5 1.5}
  \picline{4 0}{3.5 0.5}
  \picline{4 2}{3.5 1.5}
  \piccirclearc{2 0.5}{1.3}{45 135}
  \piccirclearc{2 1.5}{1.3}{-135 -45}
  \pictranslate{1 1}{
    \picrotate{90}{
      \picscale{-1 1}{
	\picfilledcircle{0 0}{0.8}{P}
      }
    }
  }
  \picfilledcircle{3 1}{0.8}{Q}
} &
\diag{6mm}{2.4}{3}{
  \piclinewidth{50}
  \picmultigraphics{2}{0 2}{
    \picellipse{1.2 0.5}{1.2 0.5}{}
  }
  \picfilledcircle{1.2 1.5}{1}{P}
} & \eepsfs[2cm]{t1-rat234}
\\[9mm]
\text{ sum $P,Q$ } & \text{ product $P\,Q$ } &
\text{ closure $\overline{P}$ } & -2\ -3\ 4\ 2
\end{array}
\]
\caption{Conway's tangles and operations with them. (The designation
`product' is very unlucky, as this operation is neither commutative,
nor associative, nor is it distributive with `sum'. Also, `sum'
is associative, but not commutative.)\label{fig1}}
\end{figure}

% \begin{figure}[htb]
% \[
% \epsfs[3cm]{t1-rat234}
% \]
% \caption{The rational tangle $4\ 3\ -2$.\label{figdgr}}
% \end{figure}
% 

\begin{defi}
A \em{tangle diagram} is a diagram consisting of strands crossing
each other, and having 4 ends. A \em{rational} tangle diagram is the
one that can be obtained from the primitive Conway tangle
diagrams by iterated left-associative product in the way displayed
in figure \reference{fig1}. (A simple but typical example of
is shown in the figure.)
\end{defi}

Let the \em{continued} (or \em{iterated}) \em{fraction}
$[[s_1,\dots,s_r]]$ for integers $s_i$ be defined inductively by
$[[s]]=s$ and
\[
[[s_1,\dots,s_{r-1},s_r]]=s_r+\frac{1}{[[s_1,\dots,s_{r-1}]]}\,.
\]
The rational tangle $T(p/q)$ is the one with Conway notation
$c_{1}\ c_{2}\ \dots c_{n}$, when the $c_i$ are chosen so that
\begin{eqn}\label{ci}
[[c_{1},c_{2},c_{3},\dots,c_n]]=\frac{p}{q}\,.
\end{eqn}
One can assume without loss of generality that $(p,q)=1$, and $0<q<|p|$.
A \em{rational (or 2-bridge) link} $S(p,q)$ is the closure of $T(p/q)$.

Montesinos links (see e.g. \cite{LickThis}) are generalizations of
pretzel and rational links and special types of arborescent links. They
are denoted in the form $M(\frac{q_1}{p_1},\dots,\frac{q_n}{p_n},e)$,
where $e,p_i,q_i$ are integers, $(p_i,q_i)=1$ and $0<|q_i|<p_i$.
Sometimes $e$ is called the \em{integer part}, and $n$ the \em{length}
of the Montesinos link. If $e=0$, it is omitted in the notation.

To visualize the link,
let $p_i/q_i$ be continued fractions of rational tangles
$c_{1,i}\dots c_{n_i,i}$ with $[[c_{1,i},c_{2,i},c_{3,i},\dots,
c_{l_i,i}]]=\frac{p_i}{q_i}$. Then $M(\frac{q_1}{p_1},
\dots,\frac{q_n}{p_n},e)$ is the link that corresponds to the
Conway notation
\[
(c_{1,1}\dots c_{l_1,1}), (c_{1,2}\dots c_{l_2,2}), \dots,
(c_{1,n}\dots c_{l_n,n}), e\,0\,.
\]
The defining convention is that all $q_i>0$ and if $p_i<0$, then
the tangle is composed so as to give a non-alternating sum with
a tangle with $p_{i\pm 1}>0$. This defines the diagram up to mirroring.

An easy exercise shows that if $q_i>0$ resp. $q_i<0$, then 
\begin{eqn}\label{Mi}
M(\dots,q_i/p_i,\dots,e)\,=\,
M(\dots,(q_i\mp p_i)/p_i,\dots,e\pm 1)\,,
\end{eqn}
i.e. both forms represent the same link (up to mirroring).

In the following the mirroring convention in the notation 
$M(q_1/p_1,\dots,q_n/p_n,e)$ will be so that a $+1$ integer twist
in $e$ is a crossing whose $A$-splicing gives an $\infty$-tangle.
So for example, the positive (right-hand) trefoil is $M(-3)$.
(If $e=0$, use \eqref{Mi} to make it $\pm 1$, and specify the
mirroring accordingly.)

A typical example is shown on figure \reference{fig2}.

\setbox\@tempboxa=\hbox{$M(\myfrac{3}{11},-\myfrac{1}{4},
\myfrac{2}{5},4)$}

\begin{figure}[htb]
\[
\begin{array}{c}
\epsfsv{5cm}{t1-dioph3} \\
\end{array}
\]
\caption{The Montesinos knot \unhbox\@tempboxa\ with
Conway notation $(213,-4,22,40)$.\label{fig2}}
\end{figure}

If the length $n<3$, an easy observation shows that the
Montesinos link is in fact a rational link.
A \em{pretzel link} is a Montesinos link with all $|q_i|=1$.
Geometric properties of Montesinos links are discussed in detail in
\cite{BurZie}.

\section{5-equivalence of Montesinos links}

\begin{defi}
We denote the equivalence of polynomials in $\bZ[t^{\pm 1/2}]$
modulo $X=(t^5+1)/(t+1)$ and multiplication with $\pm t^{k/2}$
by $\doteq$. By $\bar V$ we denote the equivalence class of a
polynomial $V$ under $\doteq$. For $c\in \bC$, we denote by
$\tl c$ the set of complex numbers obtained from $c$ by 
multiplication with a 20-th root of unity. If $\tl c_1=\tl c_2$,
we write also $c_1\simeq c_2$.
\end{defi}

\begin{prop}(Przytycki-Ishiwata \cite{Ishiwata}) 
$\bar V$ is an invariant of $5$-moves.
\end{prop}

The roots of $X=(t^5+1)/(t+1)$ are the primitive 10-th roots of unity.
These are two pairs of conjugate complex numbers. Since, when
working with real coefficients, conjugate complex numbers are
equivalent, we have two roots to work with, $e^{\pi i/5}$ and 
$e^{3\pi i/5}$. In particular we obtain from the previous
proposition:

\begin{corr}(Przytycki-Ishiwata) 
For a link $L$ and $n=1,3$, let $x_n(L)=V_L(e^{n\pi i/5})$.
Then the quantities $\gm_n=\wt{x_n}$, and in particular 
$v_n=|x_n|$ are invariants of $5$-moves.
\end{corr}

Concerning $Q$, we have the following easy to show fact (see
for example lemma \reference{lmq} below).

\begin{prop} \label{JR}
$Q(2\cos 2\pi /5)$ and $Q(2\cos 4\pi /5)$ are invariants of $5$-moves.
\end{prop}

These values of $Q$ are studied by Jones \cite{Jones2} and Rong
\cite{Rong}, who show that they are equal to $\pm \sqrt{5}^{d(K)}$,
where $d(L):=\dim_{\bZ_5}H_1(D_L,\bZ_5)$ is the number of torsion
numbers divisible by 5 of the homology of the 2-branched cover $D_L$.
This special form makes
the $Q$ values of limited use as detectors of $5$-move (in)equivalence, 
but nonetheless one should not write them off. In particular the sign
can be useful in some cases.

Ishiwata \cite{Ishiwata} succeeded in determining the $5$-move equivalence classes
of rational tangles.

\begin{theo}\label{thrt}(Ishiwata) 
Rational tangles are $5$-move equivalent to one of the 12 basic tangles
\[
0,\,\infty,\,\pm 1,\,\pm 2,\,\pm 1/2,\,\pm 3/2,\,2/5,\,5/2\,.
\]
\end{theo}

One obtains immediately from that the result for rational (2-bridge) 
links.

\begin{corr}\label{crrt}(Ishiwata) 
Rational links are $5$-move equivalent to one of the unknot, 2-component
unlink, figure-8 knot or Hopf link.
\end{corr}

This means in particular that $v_1$ and $v_3$ take only 4 values on
rational links. (They are indeed different for the 4 classes, so the
classes are distinct.) Nonetheless, some very strong conditions on 
polynomial invariants of 2-bridge links are known. For the Alexander 
polynomial such a condition is formulated in \cite{Murasugi}. The 
test of Kanenobu \cite{Kanenobu} using $V$ and $Q$ is particularly 
efficient; see also \cite{kan}. Here we use some of Ishiwata's work to 
extend her result to Montesinos links. Our goal is to obtain conditions 
on the polynomial invariants of Montesinos links (which are much
fewer in previous literature, and less efficient). 

\begin{theo}\label{thq}
A Montesinos link is equivalent up to 5-moves and mutations to 
one of
\begin{enumerate}
\item $M(-1/2,\dots,-1/2,1/2,\dots,1/2)$ with $k\ge 0$ entries 1/2
and $0\le l\le 4$ entries $-1/2$, s.t. $l+k\ge 3$,
\item\label{F2} $M(1/2,\dots,1/2,2/5,\dots,2/5)$ with $k\ge 0$
entries 1/2, and $1\le l$ entries 2/5, s.t. $l+k\ge 3$,
\item\label{F3} the connected sum of some number (possibly no) of
Hopf links, with possible trivial split components,
\item\label{F4} the connected sum of some number (at least one)
of figure-8 knots with a link of Form \reference{F3}, or
\item $M(1/2,1/2,1/2,1)$.
\end{enumerate}
\end{theo}

\proof We apply first theorem \ref{thrt} to reduce all
tangles to the 12 basic ones. If one of the tangles is
$\infty$, we are down to the connected sum of rational
links, and with corollary \ref{crrt} have the third or
fourth form.

So we assume henceforth that we have no $\infty$ tangle.
Then the 11 remaining basic tangles have non-integer parts
$\pm 1/2$ and $2/5$. We assume that at least three non-integer 
parts occur, otherwise we have again a rational link.

Now assume that at least one tangle $2/5$ occurs. Then we use the 
observation of Ishiwata (made in the proof of theorem \ref{thrt}) 
that the $2/5$ tangle acts (up to 5-moves) as an ``annihilator'' 
of integer twists. Since one can make all $-1/2$ to $1/2$ up to
integer twists, we have the second form.

So assume below that we have only $\pm 1/2$ and possible integer 
twists. It is an easy observation that 5 copies of $1/2$ are
5-move equivalent to 5 copies of $-1/2$, so we can always make 
the $-1/2$ to be at most 4. Moreover, the normal form of Montesinos 
links teaches that if have $n\ne 0$ integer twists, then we can 
assume that all $\pm 1/2$ have the sign of $n$. (If $n=0$ we have
the first form.)

Assume first $n>0$. Clearly we can achieve that $n<5$. Now if the 
number $l$ of `$1/2$' satisfies $l\ge 5-n$, we can make the $n$ 
integer twists to $5-n$ of opposite sign, and annihilate them
making $5-n$ of the $l$ copies of $1/2$ to $-1/2$. Since $0<n<5$ 
and $l\ge 3$, there is only one case where this procedure does not 
work, namely when $l=3$ and $n=1$. In this case we have the link 
of the fifth form. Otherwise we achieve the first form. 

If $n<0$ we can apply the same argument, since the family of
links of the first form is invariant (by making the $-1/2$
to be $<5$ under 5-equivalence) to the family of their mirror 
images. Also, the link in the fifth form (as well as all links
of Forms \ref{F2}, \ref{F3} and \ref{F4}) is equivalent to its
mirror image. \qed

\begin{rem}\label{rempre}
By omission of Form 2, we obtain an analogous result for pretzel
links. The (not very natural) distinction of Forms \ref{F3} and
\ref{F4} will become clear with theorem \reference{thfr} below.
\end{rem}

\section{Jones polynomial of Montesinos links}

\subsection{Formulas for the Jones polynomial}

We succeeded in evaluating $\bar V$ on the representatives
in theorem \ref{thq}, and so we have

\begin{theo}\label{thfr}
For a Montesinos link $L$, the reduced Jones polynomial 
$\bar V_L$ equals one of
\begin{enumerate}
\item Form 1, with $k+l\ge 3$, $k,l\ge 0$, $l\le 4$
\[
V_1(k,l)=\frac{(-t-1/t)^{k+l}+(-1/t)^l(1-t)^{k+l}(t+1+1/t)}
  {-\sqrt{t}-\frac{1}{\sqrt{t}}}
\]
\item Form 2, with $k+l\ge 3$, $k,l-1\ge 0$
\[
V_2(k,l)=\frac{(1-t^2)^{l}(1-t)^{k}(t+1+1/t)}{-\sqrt{t}-\frac{1}{\sqrt{t}}}
\]
\item Form 3, with $k,l\ge 0$
\[
V_3(k,l)=\left(-\sqrt{t}-\frac{1}{\sqrt{t}}\right)^l
  \left(-{t}-\frac{1}{{t}}\right)^k
\]
\item Form 4: 0
\item Form 5:
\[
V_5=\frac{(-t-1/t)^3-(1-t)^3t(t+1+1/t)}{-\sqrt{t}-\frac{1}{\sqrt{t}}}
\]
\end{enumerate}
\end{theo}

The proof of this theorem bases on routine Kauffman bracket skein
module calculation, and we omit it. (The reader may consult
\cite{gwg} for some explanation on this kind of calculation.)
The given polynomials correspond (up to units) to the Jones
polynomials of the knots, except for Form 2, where the numerator
was simplified by reducing modulo $X=(t^5+1)/(t+1)$. (For the
calculation, we also used $M(-1/2,\dots,-1/2,2/5,\dots,2/5)$,
which is equivalent.) The fraction 
may therefore not be a genuine polynomial in $\bZ[t^{\pm 1}]$. Instead 
it should be considered lying in the field $\bZ[t^{\pm 1}]/\,<X>$
(where the division is possible).

\begin{rem}\label{remhom}
We should remark the following relation of these forms to $d(L)$.
\begin{enumerate}
\item Form 1: $d(L)\le 1$. Whether $d(L)=0$ or $1$ can then be 
decided by looking at the ($5$-divisibility of the) determinant. By 
a simple calculation we find 
\[
d(L)=\,\left\{\begin{array}{cl}
0 & \mbox{if $5\nmid k-l$} \\ 1 & \mbox{if $5\mid k-l$}
\end{array}\right.
\]

\item Form 2: $|d(L)-l|\le 1$. This follows because by a
  change from a 0 tangle to an infinity tangle, $d(L)$ is
  altered by at most $\pm 1$, and we can obtain the connected 
  sum $L'$ of $k$ Hopf links and $l$ figure-8-knots, and $d(L')=l$. 
  % Now assume $k\ge 5$. Then we can also splice $0\to infty$ 
  % so that $L''=M(1/2,\dots,1/2)\# M(1/2,\dots,1/2,2/5,\dots,2/5)$
  % (with $k$ copies of
  % $1/2$ and $l$ copies of $2/5$); a calculation of the 
  % determinant of former factor shows that it is divisible by
  % 5 iff $5\mid k$. In that case, we apply the previous argument
  % on $L'=M(2/5,\dots,2/5)$ $d(L'')\ge  

\item Form 3: $l=d(L)$

\item Form 4: $d(L)>0$
\item Form 5: By direct calculation the determinant is 20, so $d(L)=1$.
\end{enumerate}
\end{rem}

For the following calculations, it is useful to compile
a few norms of complex numbers that will repeatedly occur.
\begin{eqn}\label{vals}
\begin{array}{c|c|c}
\mbox{norm of} & t=e^{\pi i/5} & t=e^{3\pi i/5} \\
\hline
\ry{1.8em}1-t &    \frac{\sqrt{5}-1}{2}\approx 0.618033988     
              & \sqrt{2+\sqrt{3-\sqrt{5}}}\approx 1.6952970385 \\[2mm]
1+t &    \frac{\sqrt{2\sqrt{5}+10}}{2}\approx 1.9021130325 
              & \sqrt{2-\sqrt{3-\sqrt{5}}}\approx 1.06111637019 \\[2mm]
1+t^2 &  \frac{\sqrt{5}+1}{2}\approx 1.618033988  
              &  \sqrt{3-\sqrt{5}} \approx 0.874032048897 \\[2mm]
1-t^2 &  \sqrt{\frac{5-\sqrt{5}}{2}}\approx 1.17557050458
              & \frac{\sqrt{2\sqrt{5}+10}}{2}\approx 1.9021130325 \\[2mm]
1+t+t^2 & \frac{\sqrt{5}+3}{2}\approx 2.618033988 
              & 1- \sqrt{3-\sqrt{5}} \approx 0.125967951102 
\end{array}
\end{eqn}

Using this table, one easily sees

\begin{corr}
The set 
\[
\{(|V(e^{\pi i/5})|,|V(e^{3\pi i/5})|)\,:\,\mbox{$L$ is a 
  Montesinos link}\,\}
\]
as a subset of $\bR^2$ is discrete.
\end{corr}

Discrete is to be understood here so that the intersection with
any ball is finite. (It is in fact in size at most logarithmic 
in the radius of the ball.)

\proof We can treat the families $V_{1,2,3}$ separately. In case of 
$V_{1,3}$, the invariant $v_1$ diverges (exponentially) when $k+l
\to \infty$. For $V_2$ take $v_3$. \qed

The pretzel links are those for which Form 2 does not occur, and so 
we have

\begin{corr}
The set 
\[
\{|V(e^{\pi i/5})|\,:\,\mbox{$L$ is a pretzel link}\,\}\,\subset\, \bR
\]
is discrete. \qed
\end{corr}

Note that this gives a very strong condition on the polynomials of
pretzel or Montesinos links. For a general link, there is no particular 
feature of $v_1$ or $v_3$ to be expected. For arbitrary links they will 
be dense in $\bR_+$. Jones claimed in \cite{Jones} (and I wrote down 
an argument in \cite{gen2}) that $v_1$ is dense in an interval on 
closed 4-braids. On the opposite hand, while several other conditions 
on invariants of Montesinos links are known, these are all difficult 
to test and/or assume additional properties of the diagram. 

For example a link $L$ is Montesinos if and only if $\pi_1(D_L)$ 
is finite (where $D_L$ is the double branched cover). While it is easy
to gain a presentation of $\pi_1(D_L)$ from any diagram of $L$, the
finiteness decision problem from this presentation is highly difficult. 

Also, the work in \cite{LickThis,Thistle} implies that any Montesinos
link has a minimal (in crossing number) Montesinos diagram. However, it
is not clear how to find this diagram starting from a given one. If the
link is alternating, all minimal diagrams will be Montesinos diagrams.
However, for non-alternating knots there may be even non-Montesinos
minimal crossing diagrams; the knots $10_{145}$, $10_{146}$
and $10_{147}$ are examples. 

In contrast, our result allows for a (not completely exact but)
efficient test of the Montesinos property. It is in fact this
condition, arising from Ishiwata's result and the detection of
5-moves by the polynomials, that motivated the present work. We
give some examples of the application of our condition in \S\ref{SA}.
First we show, however, that the invariants we have in hand suffice to
distinguish all equivalence classes in theorem \reference{thq}.

\subsection{Completing the classification}

\begin{theo}\label{thq2}
No two forms in theorem \ref{thq} are 5-move equivalent.
In other words, the links in theorem \ref{thq} determine exactly
the equivalence classes of Montesinos links under mutation and
5-moves.
\end{theo}

\begin{rem}
The massive failure of polynomial invariants to detect mutation
results in a serious lack of tools to examine 5-equivalence of 
mutants. While it is almost certain that such examples exist, probably
no simpler invariants than the (almost incomputable) Burnside
groups \cite{DP} could be applied to show it. However, note that
at least it is known that any mutant of a Montesinos link is again 
a Montesinos link, so that the study of 5-equivalence classes of
mutants among Montesinos links would reduce to the mutants of the 
representative links in theorem \ref{thfr}. This is (basically) the
question whether the 1/2,2/5 Montesinos tangle (the one with Conway
notation 2,22) is 5-equivalent to its
(mutant) 2/5,1/2 Montesinos tangle. In particular, our classification
of 5-move equivalence of pretzel links (remark \reference{rempre})
applies without mutation.
\end{rem}

\proof We check case by case possible duplications of
$\bar V_L$ and $d(L)$. Often $v_1$ and $v_3$, or sometimes
at least $\gm_1$ and $\gm_3$, already suffice to distinguish
the forms.
We will show that all these invariants coincide only if the
links are identical (except the final case, where we are 
lead to apply the $Q$ polynomial to arrive at this 
conclusion). We remark also that $V_5=V_1(4,-1)$. Although 
the value $l=-1$ in Form 1 does not make sense knot-theoretically,
from the point of view of formal algebra this identity will
allow us to handle Form 5 often in the same way as Form 1.

We also make the convention that is a comparison `Form x vs Form y',
the integers $k,l$ are parameters that correspond to `Form x',
while $k',l'$ correspond to `Form y'.

\begin{itemize}
\item Form 1 vs Form 1:
      We need to deal with 
\begin{eqn}\label{V1}
V_1(k,l)\doteq V_1(k',l')\,.
\end{eqn}
Looking at $t=e^{\pi i/5}$, we find using the above
values and $k+l\ge 3$,
\[%begin{eqn}\label{4'}
\left|\,(-1/t)^l(1-t)^{k+l}(t+1+1/t)\,\right|\,\le\,
\left(\frac{\sqrt{5}-1}{2}\right)^3\frac{\sqrt{5}+3}{2}\,=\,
\frac{\sqrt{5}-1}{2}\,.
\]%end{eqn}
Using this inequality and the analogous one with
$k,l$ replaced by $k',l'$, and taking norms in \eqref{V1}
we find
\begin{eqn}\label{est}
\left|\,\left|-t-\frac{1}{t}\right|^{k+l}-
  \left|-t-\frac{1}{t}\right|^{k'+l'}\,\right|\,\le\,\sqrt{5}-1\,.
\end{eqn}
But $\left|-t-\frac{1}{t}\right|=\frac{\sqrt{5}+1}{2}$, and
$k+l,k'+l'\ge 3$, so we see that $k+l=k'+l'$. We will argue
that then \eqref{V1} gives $5\mid l-l'$, so $(k,l)=(k',l')$, 
as desired. 

Indeed, set $\tl k=k+l=k'+l'$. Now it is more convenient to use 
$\gm_n$ instead of $v_n$. So assume there are two different values 
$l_{1,2}$ of $l=0,\dots,4$, such that $V_1({\tl k-l,l})$ (are 
different but) differ by a 20-th root of unity. Now
\[
|t+1|\,\cdot\,\big|\,V_1({\tl k-l,l})\,\big|\,\ge\,\left|\,-t
-\frac 1t\,\right|^{\tl k}-|1-t|^{\tl k}\,\left|\,1+t+\frac 1t\,\right|\,.
\]
So $V_1({\tl k-l_1,l_1})\simeq V_1({\tl k-l_2,l_2})$ implies that
\[
|t+1|\,\cdot\,\big|\,V_1({\tl k-l_1,l_1})-V_1({\tl k-l_2,l_2})
\,\big|\,\ge\,\big|\,1-e^{\pi i/10}\,\big|\,\cdot\,\left[\,
\left|\,-t-\frac 1t\,\right|^{\tl k}-
|1-t|^{\tl k}\,\left|\,1+t+\frac 1t\,\right|\,\right]\,.
\]
On the other hand, by the definition of $V_1$,
\[
|t+1|\,\cdot\,\big|\,V_1({\tl k-l_1,l_1})-V_1({\tl k-l_2,l_2})\,\big|\,\le\,
\frac{\sqrt{2\sqrt{5}+10}}{2}\,|1-t|^{\tl k}\,\left|\,1+t+\frac 1t\,\right|\,,
\]
where the first factor on the right stands for the largest distance between
two 5-th roots of unity (which are $-1/t$). Setting $t=e^{\pi i/5}$,
and combining the last two inequalities, we obtain
\[
\frac{\sqrt{2\sqrt{5}+10}}{2}\cdot \left(\,\frac{\sqrt{5}-1}{2}\right)^{\tl k}
\frac{\sqrt{5}+3}{2}\,\ge\,\big|\,1-e^{\pi i/10}\,\big|\,\cdot\,\left[\,
\left(\,\frac{\sqrt{5}+1}{2}\right)^{\tl k}-\left(\,\frac{\sqrt{5}-1}{2}\right)^{\tl k}
\,\cdot\,\frac{\sqrt{5}+3}{2}\,\right]\,.
\]

So with $\big|\,1-e^{\pi i/10}\,\big|=\sqrt{2-\sqrt{\frac{\sqrt{5}+5}{2}}}$, we have
\[
\left(\,\frac{\sqrt{5}+1}{2}\right)^{\tl k}\,\le\,\left(\,\frac{\sqrt{5}-1}{2}\right)
^{\tl k-2}\left[\,1+\frac{\sqrt{2\sqrt{5}+10}}{2\sqrt{2-\sqrt{\frac{\sqrt{5}+5}{2}}}}\right]\,.
\]
Evaluating the second factor on the right, we have
\[
\left(\,\frac{\sqrt{5}+1}{2}\right)^{2\tl k-2}\,\le\,7.0798\dots,
\]
so $\tl k\le 3$, that is, $\tl k=3$. These cases are easily checked directly.
By direct calculation, we find that $e_1$ and $e_3$ distinguish the
cases $l=0,3$ from those $l=1,2$. Then one verifies that 
\[
\left(\frac{V_1(l,3-l)\left(e^{\pi i/5}\right)}
           {V_1(3-l,l)\left(e^{\pi i/5}\right)}\right)^{20}\,\ne\, 1
\]
for $l=0,1$, and this case is finished.

\item Form 1 vs Form 2: We assume
\begin{eqn}\label{V12}
V_1(k,l)\doteq V_2(k',l')\,,
\end{eqn}
which means
\begin{eqn}\label{5'}
\left(-t-\frac{1}{t}\right)^{k+l}+\left(-\frac{1}{t}\right)^l(1-t)^{k+l}
\left(t+1+\frac{1}{t}\right)\doteq (1-t)^{\tl k}
\left(t+1+\frac{1}{t}\right)\,\left(\frac{1-t^2}{1-t}\right)^{l'}\,.
\end{eqn}
Using remark \ref{remhom}, we see that for a link $L$ fitting
into both forms we have $d(L)\le 1$, and so $l'\le 2$.
For $t=e^{\pi i/5}$, the base of the rightmost exponent in
\eqref{5'} is of norm
$>1$. We take norms in \eqref{5'} and bring the second summand
on the l.h.s. to the right. We obtain analogously to \eqref{est}
\[
\left|-t-\frac{1}{t}\right|^{k+l}\,\le\,\frac{\sqrt{5}-1}{2}\,\left(
1+\frac{2\sqrt{5}+10}{4}\,\right)\,\approx\,2.85\dots\,<\,
\left(\frac{\sqrt{5}+1}{2}\right)^3\,,
\]
which is impossible.

% \begin{itemize}
% 
% \item $l'=0$. So
% \[
% \left(-t-\frac{1}{t}\right)^{k+l}+\left(-\frac{1}{t}\right)^l(1-t)^{k+l}
% \left(t+1+\frac{1}{t}\right)\doteq (1-t)^{k}\left(t+1+\frac{1}{t}\right)\,.
% \]
% Since $k+l,k'=k'+l'\ge 3$, the norm of the term on the right
% and the second summand on the left for $t=e^{\pi i/5}$ are at most 
% $\frac{\sqrt{5}-1}{2}$, so we have again
% \[
% \left|-t-\frac{1}{t}\right|^{k+l}\,\le\,\sqrt{5}-1
% \]
% which is impossible.
% 
% \item $l'=1$. Let $\tl k=k'+l'=k'+1\ge 3$. Then we have
% \[
% \left(-t-\frac{1}{t}\right)^{k+l}+\left(-\frac{1}{t}\right)^l(1-t)^{k+l}
% \left(t+1+\frac{1}{t}\right)\doteq 
% (1-t)^{\tl k}\left(t+1+\frac{1}{t}\right)\frac{1-t^2}{1-t}\,.
% \]
% Now the norm on the right evaluates to $\frac{\sqrt{5}+3}{2}$,
% and \eqref{est} becomes 
% \[
% \left|-t-\frac{1}{t}\right|^{k+l}\,\le\,\sqrt{5}+1\,,
% \]
% which still allows the same contradiction.
% 
% \item $l'=2$. In this case by remark \ref{remhom} we must
% have $d(L)=1$, so $5\mid k-l$. Since the possible $k,l$ with 
% $k+l=3$ do not satisfy this condition, we can assume $k+l\ge 4$.
% The previous sort of calculation now leads to
% \[
%   \left|-t-\frac{1}{t}\right|^{k+l}\,\le\,
%   \frac{\sqrt{5}-1}{2}+\left(\frac{\sqrt{5}+3}{2}\right)^2
%   \frac{2}{\sqrt{5}-1}\,=\,2+2\sqrt{5}\,\approx\,6.47\,,
% \]
% while with $k+l\ge 4$ we have $\left|-t-\frac{1}{t}\right|^{k+l}
% \ge 6.85$, and so are done. 
% \end{itemize}

\item Form 1 vs Form 3: Let $\tl k=k'+l'\ge 0$. We must consider
\begin{eqn}\label{7}
  \left(-t-\frac{1}{t}\right)^{k+l}+\left(-\frac{1}{t}\right)^l(1-t)^{k+l}
  \left(t+1+\frac{1}{t}\right)\,\doteq\,
  \left(-t-\frac{1}{t}\right)^{\tl k}\left(-\sqrt{t}-\frac{1}
    {\sqrt{t}}\right)
  \left(\frac{-\sqrt{t}-\frac{1}{\sqrt{t}}}{-t-\frac{1}{t}}\right)^{l'}\,.
\end{eqn}

By remark \ref{remhom} we have $l'=0$ if $5\nmid k-l$ and
  $l'=1$ otherwise. So
\[
|1-t|^{k+l}\left|1+t+\frac{1}{t}\right|\,\le\,
|t^2+1|^{k+l}+|t+1|\cdot |t^2+1|^{\tl k}\,\max\left(\,1,\,
  \frac{|t+1|}{|t^2+1|}\,\right)\,.
\]
Look first at $t=e^{3\pi i/5}$. Since $|t+1|=\sqrt{2-\sqrt{3-\sqrt{5}}}>1$ 
  and $|t^2+1|=\sqrt{3-\sqrt{5}}<1$, we have the second alternative in 
  the maximum. Then
\[
\sqrt{2+\sqrt{3-\sqrt{5}}}^{k+l}\cdot \left(1- \sqrt{3-\sqrt{5}}\right)\,\le\,
  \sqrt{3-\sqrt{5}}^{k+l}+\sqrt{3-\sqrt{5}}^{\tl k-1}\left(2-\sqrt
    {3-\sqrt{5}}\right)\,.
\]
So
\[
\sqrt{2+\sqrt{3-\sqrt{5}}}^{k+l}\,\le\,\frac{\sqrt{3-\sqrt{5}}^{\,4}+
  \left(2-\sqrt{3-\sqrt{5}}\right)}{\left(1- 
    \sqrt{3-\sqrt{5}}\right)\,\sqrt{3-\sqrt{5}}}\,,
\]
and by calculation
\[
1.6952970385^{k+l}\,\le\, 15.5273358421\,,
\]
so $k+l\le 5$. Now look at $t=e^{\pi i/5}$. Taking norms
in \eqref{7}, and this time minimizing w.r.t. $l'=0,1$, we find
\[
\left(\frac{\sqrt{5}+1}{2}\right)^{k+l}+\left(\frac{\sqrt{5}-1}
    {2}\right)^{k+l}
  \frac{\sqrt{5}+3}{2}\,=\,
  \left(\frac{\sqrt{5}+1}{2}\right)^{k+l}+\left(\frac{\sqrt{5}+1}
    {2}\right)^{2-k-l}
  \,\ge\,\left(\frac{\sqrt{5}+1}{2}\right)^{\tl k} 
    \frac{\sqrt{2\sqrt{5}+10}}{2}\,.
\]
The biggest l.h.s. evaluates at $k+l=5$ to $\approx 11.326237921249$,
and then
\[
\left(\frac{\sqrt{5}+1}{2}\right)^{\tl k}\,\le\,\frac{2\cdot 
  11.326237921249} {\sqrt{2\sqrt{5}+10}}\,\approx\, 5.9545556584\,,
\]
so $\tl k\le 3$. So it
remains to check the cases $k'+l'\le 3$, $3\le k+l\le 5$ (with $k,l,k'
\ge 0$, $l\le 4$; and $l'=1$ if $5\mid k-l$ and $l'=0$ otherwise). 
It is easy to perform (still better by computer) these handful of
comparisons; $v_3$ distinguishes all such $V_1(k,l)$ and $V_3(k',l')$.

\item Form 1 vs Form 4: The same estimate as with 
  'Form 1 vs Form 1' works. 

\item Form 1 vs Form 5: We can apply the same argument as with 
'Form 1 vs Form 1', since we did not use that $l>0$ there. 
Then again we need to check only $\tl k=3$. Direct calculation
shows that $e_1$ and $e_3$ distinguish the case $l=-1$ (in Form 5) 
from $l=0,3$ and $l=1,2$. 
% Then one verifies that 
% \[
% \left(\frac{V_1(l,3-l)\left(e^{\pi i/5}\right)}
%            {V_1(3-l,l)\left(e^{\pi i/5}\right)}\right)^{20}\,\ne\, 1
% \]
% for $l=0,1$, and this case is finished.
% 

\item Form 2 vs Form 2: 
\[
(1-t^2)^l(1-t)^k\,\doteq\,(1-t^2)^{l'}(1-t)^{k'}\,.
\]
Since $|l-d(L)|,|l'-d(L)|\le 1$, we have $|l-l'|\le 2$. Since for
$l=l'$ we are easily done, assume w.l.o.g. $l'-l\in\{1,2\}$.
Then
\[
|1-t|^{k-k'}\in \left\{\,\frac{1}{|1-t^2|},\, \frac{1}{|1-t^2|^2}\,
 \right\}\,.
\]
For $t=e^{\pi i/5}$ the numbers on the right are $0.850\dots$, $0.723\dots$,
while the sequence on the left is $1$, $0.618\dots$, $0.381\dots$, 
$0.236\dots$ etc.

\item Form 2 vs Form 3: We have to check
\begin{eqn}\label{uo}
(1-t^2)^{l}(1-t)^{k}\left(t+1+\frac{1}{t}\right)\,\doteq\,
\left(-\sqrt{t}-\frac{1}{\sqrt{t}}\right)^{l'+1}
  \left(-{t}-\frac{1}{{t}}\right)^{k'}\,.
\end{eqn}
Now by remark \ref{remhom}, $|l-l'|\le 1$, so $\tl l:=l'+1-l=0,1,2$.
Taking norms in \eqref{uo} and using $\tl l=0,1,2$, we have
\[
|1-t|^{l+k}\,=\,\frac{|t+1|^{\tl l}|t^2+1|^{k'}}
  {\ds\left|1+t+\frac{1}{t}\right|}\,.
\]
Now for $t=e^{\pi i/5}$, we get
\[
\left(\frac{\sqrt{5}-1}{2}\right)^{k+l}\,=\,\frac
  {\left(\ds\frac{\sqrt{5}+1}{2}\right)^{k'}\cdot
   \left(\frac{\sqrt{2\sqrt{5}+10}}{2}\right)^{\tl l}
  }
  {\ds \frac{\sqrt{5}+3}{2}}\,=\,
  \left(\frac{\sqrt{5}+1}{2}\right)^{k'-2}\cdot
   \left(\frac{\sqrt{2\sqrt{5}+10}}{2}\right)^{\tl l}\,.
\]
Now the r.h.s. is at least $\left(\frac{\sqrt{5}-1}{2}\right)^2$,
  while the l.h.s. for $k+l\ge 3$ is at most $\left(\frac{\sqrt{5}
  -1}{2}\right)^3$, a contradiction.

% ***********************************************************

\item Form 2 vs Form 4: trivial

\item Form 5 vs Form 2: We apply the same argument as with `Form 1 vs Form 2'.
%   This works except in the case $l'=2$, where we use that $k+l\ge 4$. (Now
%   indeed $5\mid k-l$.) Assuming $l'=2$, now look in
% \[
% \left(-t-\frac{1}{t}\right)^{3}-t(1-t)^3\left(t+1+\frac{1}{t}\right)\doteq 
% (1-t)^{\tl k}\left(t+1+\frac{1}{t}\right)\left(\frac{1-t^2}{1-t}\right)^2\,.
% \]
% at $t=e^{3\pi i/5}$. We have
% \[
% \left|t+1+\frac{1}{t}\right|\,\cdot\,|1-t|^{k'+l'}\,\le\,
% |t^2+1|^{k'+l'}+|t^2+1|^{3}+|1-t|^3\left|t+1+\frac{1}{t}\right|\,.
% \]
% Now $k'+l'\ge 3$. So we have
% \[
% 1.6952970385^{k'+l'}\,\le\,\frac{2(\sqrt{3-\sqrt{5}})^3}{1-\sqrt{3-\sqrt{5}}}+
%   \left(\sqrt{2+\sqrt{3-\sqrt{5}}}\right)^3 \approx 15.47\,,
% \]
% whence $k'+l'\le 5$. Since $l'=2$ and $3\le l'+k'$, there
% remain to check $l'=2$ and $1\le k'\le 3$. Again a simple
% computer calculation shows that $v_3$ differs on $V_5$ and
% any of the three possible $V_2(k',l')$.

\item Form 3 vs Form 3: If $V_3(k,l)\doteq V_3(k',l')$, we can assume 
$|l-l'|\le 2$ as in remark \ref{remhom}, so, as in `Form 2 vs Form 2', we have
\[
\left|-\frac{1}{t}-t\right|^{k-k'}\in \left\{\,\frac{1}{|1+t|},\, 
 \frac{1}{|1+t|^2}\, \right\}\,.
\]
For $t=e^{\pi i/5}$, we have $\left|-\frac{1}{t}-t\right|\approx 
1.618033988$ and $|1+t|\approx 1.9021130325$, so the terms on the
right are approx. $0.525$ and $0.276$, while the left side ranges
in $0.618$, $0.381$, $0.236$ etc.

\item Form 3 vs Form 4: trivial 

\item Form 5 vs Form 3: By the same norm estimate as in 
`Form 1 vs Form 3', we are left with $k'+l'\le 3$, and since 
$d(L)=1$ in Form 5, also $l'=1$, so $k'\le 2$. So it
remains to check $0\le k'\le 2$, $l'=1$. A check of $v_3$
of $V_5$ and the three possible $V_3(k',l')$ rules out
$k'=1,2$. The case $k'=0$, however, is not ruled out; 
neither it is by $v_1$, and in fact not even by $\bar V$.
It is the question of 5-equivalence of the link $7^3_1$
in \cite[appendix]{Rolfsen} (which is the last link in
theorem \ref{thfr}) and the 2-component unlink. This problem
was encountered also in Ishiwata's tabulation \cite{Ishiwata}
of 5-equivalence of links up to 9 crossings. The problem is now
resolved using the $Q$ polynomial. Clearly for both links $d(L)=1$.
However, exactly the sign, which $Q(2\cos 2\pi /5)$ contains 
additionally, manages to distinguish the two links. 

\item In Form 4 vs Form 4, Form 4 vs Form 5 and Form 5 vs Form 5
there is nothing to do, and so our proof is complete. \qed 
\end{itemize}

\section{Applications\label{SA}}

The preceding discussion explains how to proceed to test the
Montesinos property of some link $L$ using the Jones polynomial.

For Form 1, one should take $t=e^{\pi i/5}$. Then $|t^2+1|>1>|1-t|$.
Increase $m=k+l$ from 3 on, as long as
\[
|t^2+1|^{m}\le |t+1|\cdot v_1(L)+|1-t|^3\cdot |t+1+1/t|\,.
\]
If 
\[
\big|\,|t^2+1|^{m}-|t+1|\cdot v_1(L)\,\big|\,\le\, |1-t|^m\cdot |t+1+1/t|\,,
\]
then compare $\bar V$ with $\bar V_1(k,l)=\bar V_1(m-l,l)$ for $0\le l\le 4$.
Note that we have a restriction on $l$ from $d(L)$. While we cannot determine 
$d(L)$ from the Jones polynomial, it still leaves a ``trace'' in form of the
condition $5^{d(L)}\mid \dt(L)=|V_L(-1)|$. So among the 5 possible $l$
we can exclude the cases where $5$ divides exactly one of $k-l=m-2l$ 
and $V_L(-1)$, but not the other.

For Form 2, one should take $t=e^{3\pi i/5}$. Then $|t^2-1|>|1-t|>1$.
Now we increase first $l$ from 1 onward as long as $5^{l-1}\mid V_L(-1)$
(because of the relation to $d(L)$), and 
\[
|t^2-1|^{l}|t+1+1/t|\,\le\,|t+1|\cdot v_3(L)\,.
\]
For such $l$, iterate $k$ from $\max(0,3-l)$ onward, as long as
\[
|t^2-1|^{l}|1-t|^k|t+1+1/t|\,\le\,|t+1|\cdot v_3(L)\,.
\]
If equality holds, compare $\bar V$ with $\bar V_2(k,l)$.

With a similar procedure one tests form 3, now using $t=e^{\pi i/5}$. 
(Then $|t^2+1|,|1+t|>1$, and we must have $5^l\mid \dt(L)$.)
Forms 4 and 5 (latter actually being redundant, since equivalent,
as far as $\bar V$ can tell, to Form 3 for $k=0,l=1$) are tested
by direct comparison.

If one fails to find $\bar V$ in these forms, one can conclude that
$L$ is not (even 5-equivalent to) a Montesinos link. I wrote a
computer program that performs this test, and show its output on
(say) non-alternating 10 crossing knots\footnote{In this table,
we use the numbering of \cite[appendix]{Rolfsen}, but the mirroring
convention in \cite{KnotScape}; a conversion to Rolfsen's mirroring
can be found on my website.}:

\begin{verbatim}
10 124 match form 3 for k=2, l=0    10 145 match form 2 for k=2, l=1
10 125 match form 3 for k=2, l=0    10 146 match form 2 for k=2, l=1
10 126 match form 3 for k=2, l=0    10 147 match form 2 for k=2, l=1
10 127 match form 3 for k=2, l=0    10 148 match form 1 for k=1, l=2
10 128 match form 1 for k=2, l=1    10 149 match form 1 for k=1, l=2
10 129 form 5                       10 150
10 130 match form 1 for k=0, l=3    10 151
10 131 match form 1 for k=1, l=2    10 152
10 132 form 5                       10 153
10 133 match form 1 for k=2, l=1    10 154
10 134 match form 1 for k=3, l=0    10 155 match form 2 for k=1, l=2
10 135 match form 1 for k=3, l=0    10 156 match form 2 for k=1, l=2
10 136 match form 2 for k=1, l=2    10 157
10 137 match form 2 for k=1, l=2    10 158 match form 2 for k=1, l=2
10 138 match form 2 for k=1, l=2    10 159
10 139 match form 3 for k=1, l=0    10 160
10 140 match form 3 for k=0, l=0    10 161 match form 2 for k=1, l=2
10 141 match form 3 for k=0, l=0    10 162 form 5
10 142 form 5                       10 163 match form 2 for k=3, l=1
10 143 match form 3 for k=1, l=0    10 164 match form 2 for k=1, l=2
10 144 match form 3 for k=0, l=0    10 165
\end{verbatim}

Whenever no form is found, the Montesinos property is ruled out. This 
happens for 9 of the last 18 knots; the first 24 knots are Montesinos.
Our test thus seems relatively efficient. It can surely not be perfect, 
since it is invariant under 5-moves, and also sporadic duplications of 
Jones polynomials occur. However, it seems easier to perform than all 
previously known Montesinos link tests (at least for general 
diagrams, and as long as the Jones polynomial can be calculated).

\section{Invariants of the Kauffman polynomial}

Now we study the 5-move invariants of the Kauffman polynomial to
enhance the test. The following is easy to see:

\begin{lemma}\label{lmq}
Let $n,m\in \bZ_k$, so that $n\ne \pm m$ and $w=e^{2\pi i n/k}
\ne \pm 1, \pm i$. Then $F(a,z)$ for $a=e^{2\pi i m/k}$ and
$z=w+w^{-1} = 2\cos 2\pi n/k$, up to multiplication with
(powers of) $a$, is a $k$-move invariant.
\end{lemma}

\proof[(sketch)]
Consider the generating function
\[
f(a,z,x)=\sum_{j=0}^\infty \Lm(A_j)(a,z)x^j\,,
\]
where $A_j$ are link diagrams with a twist tangle of $j$ crossings
(that is, a tangle with Conway notation $j$, as on the right of
\eqref{3m} for $j=3$). Use the relations
\eqref{(2)} to rewrite $f$ as a rational function in $x$,
determined by $A_{0,1,\infty}$. Finally analyze for what
values of $a$ and $z$ (for which $F(a,z)$ makes sense), the zeros
$x$ of the denominator polynomial are distinct $k$-th roots of
unity. \qed

In the case of $k=5$ this gives up to conjugacy four values.
For $m=1$ we have the two (equivalent under the interchange of
$\pm \sqrt{5}$) Jones-Rong values of $Q$ in proposition \ref{JR}.
The other two values are $F(e^{\pm 2\pi i/5},2\cos 4\pi /5)$ and
$F(e^{\pm 4\pi i/5},2\cos 2\pi /5)$. They are equivalent in the
same sense, but (as with $v_1$ and $v_3$) the behaviour of the
complex norm under this equivalence could be drastic, so one should
not \em{a priori} discard one. (One could do so, though, for one
Jones-Rong value, because its special form implies that the
$\pm \sqrt{5}$ equivalence results at most in a change of sign.)

To evaluate $F$ on the links in theorem \reference{thq}, we
make again a skein module calculation, this time using the
Kauffman polynomial (\em{not Kauffman bracket}!) skein relation \eqref{(2)}.
In the following we assume that $a$ and $z$ are as specified in
lemma \reference{lmq} for $k=5$, and consider $F(a,z)$ up to powers
of $a$.

We start with determining the coefficients
\begin{eqnarray*}
<1/2> & = & A\,<0> +  B\,<1> +  C\,<\infty> \\
<2/5> & = & \Mvariable{D}\,<0> +  \Mvariable{E}\,<1> +  \Mvariable{F}\,<\infty>
\end{eqnarray*}
of the tangles $T=<1/2>$ and $T'=<2/5>$ in the Kauffman skein module
(not to be confused with the Kauffman bracket skein module).

By taking the closure of the sum of the tangle $T=<1/2>$ resp. 
$T'=<2/5>$ with the $0$, $\infty$ and $-1$ tangles, we obtain a
linear equation system that determines the Kauffman skein module
coefficients of $T$ and $T'$.

Let first ${a_1} = \frac{1}{a} + a$, and
\begin{eqnarray*}
\Mvariable{T_2}  & = &  -1 + \frac{{a_1}}{z} \\
\Mvariable{H}    & = &  z a_1-\Mvariable{T_2} \\
\Mvariable{G_4}  & = &  (1-a_1^2)-z a_1+z^2 a_1^2+z^3 a_1 \\
\Mvariable{G_3}  & = &  (-1/a-2a)+z a_1/a+z^2 a_1
\end{eqnarray*}
be the writhe-unnormalized polynomials of the 2-component unlink,
Hopf link, figure-8 knot, and negative trefoil resp.

Let
\[
M= \left( 
  \matrix{ \Mvariable{T_2} & a & 1 \cr 1 & \frac{1}{a} & \Mvariable{T_2}
  \cr \frac{1}{a} & \Mvariable{T_2} & a \cr  }
  \right)
\]
be the matrix that represents the closures of the three
skein module generating tangles $<0>$, $<1>$ and $<\infty>$, and
\[
\Mvariable{v} \,= \, \left( 
  \matrix{ a^2 & a^2\,\Mvariable{H} \cr \Mvariable{H} &
  \Mvariable{G_4} \cr \frac{1}{a} & \Mvariable{G_3} \cr  }
  \right)
\]
be the to-result polynomials for the closures of the the
tangles $T$ and $T'$. Then
\[
 \left( \matrix{ A & B & C \cr \Mvariable{D} & \Mvariable{E}
    & \Mvariable{F} \cr  } \right)
    = \left (\, M^{-1}\,\cdot\,\Mvariable{v}\,\right)^T\,,
\]
and by calculation we find
\begin{eqnarray*}
(A,B,C) & = & (za,\,z,\,-1) \\
(\Mvariable{D},\Mvariable{E},\Mvariable{F}) & = &
  \left( -a^2 + z^2 + a^2\,z^2 + a\,z^3,\quad -z + \frac{z^2}{a} + z^3,
     \quad -\frac{z}{a}  - z^2\,\right) \,.
\end{eqnarray*}

A routine (and so omitted) further calculation leads to the formulas
we wish (the link in Form 5 can be evaluated directly).

\begin{theorem}\label{thfs}
For $a$ and $z$ as in lemma \reference{lmq} for $k=5$,
the values $F(a,z)$ of the links in theorem 
\reference{thq} are given, up to powers of $a$, as follows.

For Form 1 we have with $k=n$, $l=n'$:
\begin{eqnarray*}
F_1(n,n')& = &
  \frac{1}{\Mvariable{T_2}}\,\left[\,{\Mvariable{H}}^{n + {n'}} + 
     a^{n - {n'}}\,z^{n + {n'}}\,
      \left( \sum_{j = 0}^{n}\sum_{{j'} = 0}^{{n'}}
           {n\choose j}\,{{{n'}}\choose {{j'}}}\,
            \Muserfunction{W_1}({{\left( j - {j'} \right) }\bmod 5})
        \right)\,\right]\,,
\end{eqnarray*}
where
\begin{eqnarray*}
W_1(0)     & = & -1         + {\Mvariable{T_2}}^2\,, \\
W_1(\pm 1) & = & -a^{\mp 2} + \Mvariable{T_2}\,, \\
W_1(\pm 2) & = & -a^{\pm 1} + a^{\mp 2}\Mvariable{H}\,\Mvariable{T_2}\,.
\end{eqnarray*}

For Form 2 we have
\begin{eqnarray*}
F_2(n,n') & = &
  \frac{1}{\Mvariable{T_2}}\,\left[\,{\Mvariable{G_4}}^{{n'}}\,{\Mvariable{H}}^n + 
     \sum_{j = 0}^{n}\sum_{{j'} = 0}^{{n'}}
         {n\choose j}\,{{{n'}}\choose {{j'}}}\,A^{n - j}\,
          B^j\,{\Mvariable{D}}^{{n'} - {j'}}\,
          {\Mvariable{E}}^{{j'}}\,
          \Muserfunction{W_2}({{\left( j + {j'} \right) 
     }\bmod 5}) \right] \,,
\end{eqnarray*}
where
\begin{eqnarray*}
W_2(0)     & = & -1 + {\Mvariable{T_2}}^2\,, \\
W_2(\pm 1) & = & -a^{\mp 1} + a^{\pm 1} \,\Mvariable{T_2}\,, \\
W_2(\pm 2) & = & -a^{\mp 2} + \Mvariable{H}\,\Mvariable{T_2}\,.
\end{eqnarray*}

For Form 3 and Form 4,
\begin{eqnarray*}
F_3(n,n',n''):=
  {\Mvariable{G_4}}^n\,{\Mvariable{H}}^{{n'}}\,
   {\Mvariable{T_2}}^{{n''}}\,.
\end{eqnarray*}
If $n>0$, we have Form 4; if $n=0$, Form 3 with $k=n'$,
  $l=n''$.
\end{theorem}

These expressions look slightly more complicated than those
for $V$, but are still straightforward to evaluate.
Using the formulas, we can for example rule out $10_{155}$
and $10_{161}$ from being Montesinos. The calculation in
\S\reference{SA} shows that if $10_{155}$ were Montesinos,
it must be in the $5$-equivalence class of a Montesinos
link corresponding to the polynomial $F_2(1,2)$ (for example
$10_{136}$). But the $F$ polynomial invariants of $F(10_{155})$
are different from those of $F_2(1,2)$. The same argument rules
out $10_{161}$. The other 7 undecided knots remain, and it is
suspectable that they are $5$-equivalent to Montesinos links.

\noindent{\bf Acknowledgement.}
Main motivation for this work provided Ishiwata's talk \cite{Ishiwata}.
The calculations were partly assisted by MATHEMATICA\TM\ \cite{Wolfram},
the table access program KnotScape \cite{KnotScape} and the graphic
web interface Knotilus \cite{FR}. This research was supported by
Postdoc grant P04300 of the Japan Society for the Promotion of Science
(JSPS) at the University of Tokyo. I wish to thank also to my host
T.~Kohno for his encouragement.

\end{document}